\documentclass{amsart}

\usepackage[utf8]{inputenc} 
\usepackage[english]{babel}

\theoremstyle{plain}
\newtheorem*{theorem}{Theorem} 
\newtheorem{thm}{Theorem} 
\newtheorem{lem}[thm]{Lemma}

\theoremstyle{definition}

\theoremstyle{remark} 
\newtheorem*{rem}{Remark}

\begin{document} 
\title{On the Sidon Constant for Dirichlet Polynomials} 
\author{Ole Fredrik Brevig} 
\date{\today} 
\address{Department of Mathematical Sciences, Norwegian University of Science and Technology (NTNU), NO-7491 Trondheim, Norway} 
\email{ole.brevig@math.ntnu.no}

\begin{abstract}
	We estimate the error term in the asymptotic formula of the Sidon constant for (ordinary) Dirichlet polynomials by providing explicit lower and upper bounds. The lower bound is already implicitly known, but we supply the necessary computations to make it explicit. The upper bound is improved by a factor of $\log\log\log{x}/\sqrt{\log\log{x}}$. 
\end{abstract}

\maketitle

\section{Introduction} Given a Dirichlet series $f(s) = \sum_{n=1}^\infty a_n/n^s$ we consider the truncated Dirichlet polynomials 
\begin{equation}
	\label{eq:diripoly} f_x(s) = \sum_{n \leq x} \frac{a_n}{n^s}. 
\end{equation}
The Sidon constant for Dirichlet polynomials is defined as
\[S(x) = \sup_{f_x \not \equiv 0} \frac{\|\widehat{f_x}\|_1}{\|f_x\|_\infty},\]
where $\|\widehat{f_x}\|_1 = \sum_{n\leq x} |a_n|$ and $\|f_x\|_\infty = \sup_{t \in \mathbb{R}} |f_x(it)|$. After investigations by Queff{é}lec \cite{Queffelec}, Konyagin--Queff{é}lec \cite{KQ01} and de la Bret{è}che \cite{Br08}, the formula 
\begin{equation}
	\label{eq:sidonformula} S(x) = \sqrt{x} \exp\left(\left(-1+o(1)\right)\sqrt{\frac{\log{x}\log\log{x}}{2}}\right), 
\end{equation}
as $x \to \infty$, was finally obtained by Defant--Frerick--Ortega-Cerd{\`a}--Ouna{\"{\i}}es--Seip \cite{bh-hyper}, using their hypercontractive Bohnenblust--Hille inequality for homogenous polynomials. 

The fact that a hypercontractive Bohnenblust--Hille inequality for homogenous polynomials was the final ingredient in the proof of \eqref{eq:sidonformula} should not come as a surprise, given the history of this inequality. H. Bohr \cite{bohr-darstellung,B13a} studied the following convergence abscissas for Dirichlet series: 
\begin{align*}
	\sigma_c &= \inf \left\{\sigma \, : \, \sum_{n=1}^\infty a_n / n^{\sigma} \text{ converges}\right\} & \text{(Simple)} \\
	\sigma_b &= \inf \left\{\sigma \, : \, \sum_{n=1}^\infty a_n / n^{\sigma+it} \text{ converges uniformly for } t\in\mathbb{R}\right\} &\text{(Uniform)} \\
	\sigma_a &= \inf \left\{\sigma \, : \, \sum_{n=1}^\infty |a_n| / n^{\sigma} \text{ converges}\right\} &\text{(Absolute)} 
\end{align*}

It is clear that $\sigma_a-\sigma_c\leq 1$, and this is easily seen to be optimal by the example $a_n = (-1)^{n-1}$. Bohr also computed $\sigma_a-\sigma_b \leq 1/2$, and asked whether this was optimal. His question remained open for over a decade, until Bohnenblust--Hille \cite{BH31} established their inequality and used it to give a positive answer.

Now, under the assumption that $\sum_{n=1}^\infty a_n$ diverges, the Cahen--Bohr formulas for these abscissas are
\[\sigma_c = \limsup_{x \to \infty} \frac{\log|f_x(0)|}{\log{x}},\,\,\, \sigma_b = \limsup_{x \to \infty} \frac{\log\|f_x\|_\infty}{\log{x}}\, \text{ and } \,\sigma_a = \limsup_{x \to \infty} \frac{\log\|\widehat{f_x}\|_1}{\log{x}}.\]
Clearly, $S(x)$ is connected to the relationship between uniform and absolute convergence for Dirichlet series. In fact, the optimality of $\sigma_a-\sigma_b\leq1/2$ follows directly from \eqref{eq:sidonformula} in view of the relevant Cahen--Bohr formulas. Hence \eqref{eq:sidonformula} can be considered a stronger version of Bohnenblust--Hille's result on the optimality of $1/2$. 

The quantity $S(x)$ is sometimes called the Sidon constant for the index set $\Lambda_x = \{\log{n}\,:\, n \leq x\}$. In fact, for any finite index set of real numbers $\Lambda$, let us define the Sidon constant
\[S(\Lambda) = \inf\left\{C \, :\, \sum_{\lambda \in \Lambda} |a_\lambda| \leq C \sup_{t \in \mathbb{R}} \left|\sum_{\lambda\in\Lambda} a_\lambda e^{-i\lambda t}\right|, \, \forall a_\lambda\in\mathbb{C}\right\}.\]

The most studied index set is perhaps $\Lambda_N = \{0,\,1,\,\ldots,\,N\}$, dating at least back to Erd\H{o}s \cite{erdos}. This corresponds to the study of trigonometric polynomials of the form $P(z) = \sum_{n=0}^N a_n z^n$, where $z=e^{-it}\in \mathbb{T}$. Kahane \cite{kahane1980polynomes} obtained the asymptotic formula $S(\Lambda_N) = (1-o(1))\sqrt{N}$, which was later sharpened by Bombieri--Bourgain \cite{bombieri2009kahane} to
\[\sqrt{N}\left(1 - \mathcal{O}\left(N^{-1/9+\epsilon}\right)\right) \leq S(\Lambda_N) \leq \sqrt{N}.\]

Previously, establishing \eqref{eq:sidonformula} has been the main goal, and \cite{Br08,bh-hyper} does not provide explicit estimates for the $o(1)$-term. However, by going through the proofs in \cite{Br08,bh-hyper} carefully, one can obtain explicit bounds. The goal of this paper is to sharpen \eqref{eq:sidonformula} by providing explicit lower and upper bounds for the $o(1)$-term.
\begin{theorem}
	Let $\delta(x)$ denote the $o(1)$-term in \eqref{eq:sidonformula}. Then
	\[-\frac{1}{2}\frac{\log\log\log{x}}{\log\log{x}}+\mathcal{O}\left(\frac{1}{\log\log{x}}\right) \leq \delta(x) \leq 3\,\frac{\log\log\log{x}}{\log\log{x}}+\mathcal{O}\left(\frac{1}{\log\log{x}}\right).\]
\end{theorem}

The lower bound is implicit in \cite{Br08}, but our upper bound is improved by a factor of $\log\log\log{x}/\sqrt{\log\log{x}}$ from what is obtained by combining \cite{Br08} and \cite{bh-hyper}. Our main effort will be directed at the improved upper bound. Let us first provide the necessary computations to obtain the explicit lower bound. 
\begin{proof}
	[Proof of the lower bound] In section 2.1 of \cite{Br08}, the penultimate estimate is 
	\begin{equation}
		\label{eq:sidonlower} S(x) \gg \sqrt{\frac{x\log{y}}{\log\log{x}}}\exp\left(\frac{1}{2}\log{\rho(u)}-\frac{\log{x}}{2u}\right), 
	\end{equation}
	where $u = \log{x}/\log{y}$. Dickman's function $\rho(u)$ has expansion
	\[\rho(u) = \exp\left(-u\left(\log{u}+\log\log{u}+\mathcal{O}(1)\right)\right),\]
	and we may choose $y$ in the range $x \geq y \geq \exp\left((\log\log{x})^{3/5+\epsilon}\right)$, see (2$\cdot$2) of \cite{Br08}. The parameter $y = \exp\left(\sqrt{\log{x}\log\log{x}/2}\right)$ is chosen. This yields 
	\begin{align*}
		\log{u} &= \frac{1}{2}\left(\log\log{x}+\log{2}-\log\log\log{x}\right), \\
		\log\log{u} &= -\log{2}+\log\log\log{x}+\log{\left(1+\frac{\log{2}}{\log\log{x}}-\frac{\log\log\log{x}}{\log\log{x}}\right)}. 
	\end{align*}
	In particular $\log{u}+\log\log{u} = \left(\log\log{x}+\log\log\log{x}\right)/2+\mathcal{O}(1)$. Inserting this into \eqref{eq:sidonlower} and using that $\log{y}\geq \log\log{x}$ we have 
	\begin{align*}
		S(x) &\geq \sqrt{x} \exp\left(-\sqrt{\frac{\log{x}}{2\log\log{x}}}\left(\frac{\log\log{x}}{2}+\frac{\log\log\log{x}}{2}+\mathcal{O}(1)\right)-\frac{\log{x}}{2u}\right) \\
		&= \sqrt{x}\exp\left(\left(-1-\frac{1}{2}\frac{\log\log\log{x}}{\log\log{x}}+\mathcal{O}\left(\frac{1}{\log\log{x}}\right)\right)\sqrt{\frac{\log{x}\log\log{x}}{2}}\right) 
	\end{align*}
	as required. 
\end{proof}

\section{Preliminaries} Fix some $x \geq 2$, and let $f(s)$ denote the Dirichlet polynomial \eqref{eq:diripoly}. Bohr \cite{B13b} discovered that one can lift $f(s)$ to the polydisk $\mathbb{D}^k$, where $k=\pi(x)$. For $n\leq x$, the fundamental theorem of arithmetic allows the factorization 
\begin{equation}
	\label{eq:pfact} n = \prod_{j=1}^k p_j^{\alpha_j}. 
\end{equation}
Let $\alpha(n)=(\alpha_1,\alpha_2,\ldots,\alpha_k)$. The Bohr lift of $f(s)$ is the polynomial
\[F(z) = F(z_1,\,z_2,\,\ldots,\,z_k) = \sum_{n\leq x} a_n z^{\alpha(n)}.\]
By Kronecker's Theorem, $\|F\|_\infty = \sup_{z \in \mathbb{T}^k} |F(z)|= \|f\|_\infty$, since the logarithm of the prime numbers are rationally independent. A key insight of the Bohr lift is that each prime number corresponds to an independent variable, $p_j^{-it} \longleftrightarrow z_j$. 

To control the number of variables, let $P^+(n)$ and $P^-(n)$ denote the largest and smallest prime factor of $n$ respectively. We say that $n$ is \emph{$y$-smooth} if $P^+(n)\leq y$ and \emph{$y$-rough} if $P^-(n)>y$. The following lemma is due to Konyagin--Queff{é}lec \cite[p. 171]{KQ01}, but we include the proof for the reader's convenience.
\begin{lem}
	\label{lem:nsplit} Fix some $x\geq y \geq 2$ and write $n=\kappa\eta$ where $\kappa$ is $y$-smooth and $\eta$ is $y$-rough. Then $\|f_\kappa\|_\infty \leq \|f\|_\infty$, where 
	\begin{equation}
		\label{eq:nsplit} f(s) = \sum_{n\leq x} \frac{a_n}{n^s} = \sum_{\kappa} \left(\sum_{\eta} \frac{a_n}{\eta^s}\right)\frac{1}{\kappa^s} = \sum_{\kappa} \frac{f_\kappa(s)}{\kappa^s}. 
	\end{equation}
	\begin{proof}
		Write $w=(w_1,w_2) = \left((z_1,\ldots,z_{\pi(y)}),\,(z_{\pi(y)+1},\ldots,z_{\pi(x)})\right)$ and decompose $F(w) = F(w_1,w_2) = \sum_{\kappa} F_\kappa(w_2) w_1^{\alpha(\kappa)}$. By orthogonality,
		\[F_\kappa(w_2) = \int_{\mathbb{T}^{\pi(y)}} F(w_1,w_2) \overline{w_1}^{\alpha(\kappa)}\, d\mu^{\pi(y)}(w_1).\]
		Thus clearly $\left|F_\kappa(w_2)\right| \leq \sup_{w_1} \left|F(w_1,w_2)\right|$ and hence $\|F_\kappa\|_\infty \leq \|F\|_\infty$. 
	\end{proof}
\end{lem}

The main theorem of \cite{bh-hyper} is the hypercontractive Bohnenblust--Hille inequality for homogenous polynomials: Let $m$ and $n$ be positive integers larger than $1$. Then 
\begin{equation}
	\label{eq:bhpoly} \left(\sum_{|\alpha|=m} \left|a_\alpha\right|^\frac{2m}{m+1}\right)^\frac{m+1}{2m} \leq \left(1+\frac{1}{m-1}\right)^{m-1} \sqrt{m} \left(\sqrt{2}\right)^{m-1} \sup_{z \in \mathbb{D}^n} \left|\sum_{|\alpha|=m} a_\alpha z^\alpha\right| 
\end{equation}
for every $m$-homogenous polynomial $\sum_{|\alpha|=m} a_\alpha z^\alpha$ on $\mathbb{C}^n$.\footnote{In a recent paper \cite{bpsbohr}, appearing after the completion of the present work, it was shown that the constant in the polynomial Bohnenblust--Hille inequality is subexponential. However, this interesting fact will not improve our results.}
\begin{rem}
	In view of the homogenization procedure
	\[P(z_0,\,z_1,\,\ldots,\,z_n) = z_0^{m} Q\left(\frac{z_1}{z_0},\,\frac{z_2}{z_0},\,\ldots,\,\frac{z_n}{z_0}\right)\]
	and the maximum modulus principle, it is clear that \eqref{eq:bhpoly} in fact holds for any $m$th degree complex polynomial on $\mathbb{C}^n$. We shall not need this fact. 
\end{rem}

We want to use a version of $\eqref{eq:bhpoly}$ for Dirichlet polynomials, through the Bohr lift. Let $\Omega(n)$ denote the number of prime divisors of $n$, counting multiplicity. We have $\Omega(n) = |\alpha(n)| = \alpha_1+\alpha_2+\cdots+\alpha_k$, in view of the factorization \eqref{eq:pfact}. The following lemma is a sharper version of Theorem III-1 in \cite{Queffelec}. 
\begin{lem}
	\label{lem:homoineq} For any Dirichlet polynomial $f(s) = \sum_{n\leq x} a_n / n^s$ we have 
	\begin{equation}
		\label{eq:bh} \left(\sum_{\Omega(n) = m} |a_n|^\frac{2m}{m+1}\right)^\frac{m+1}{2m} \leq e^m \|f\|_\infty. 
	\end{equation}
	\begin{proof}
		Split the polynomial into its homogenous parts, $F(z) = \sum_{m} F_m(z)$. Fix $z\in\mathbb{T}^k$, and write
		\[F\left(e^{i\theta}z_1,\,e^{i\theta}z_2,\,\ldots,\, e^{i\theta}z_k\right) = \sum_{m} \left(e^{i\theta}\right)^m \sum_{|\alpha|=m} a_\alpha z^\alpha = \sum_{m} e^{im\theta}F_m(z).\]
		Since $F_m(z)$ appear as the Fourier coefficients, $\left|F_m(z)\right| \leq \|F\|_\infty$, and thus $\|F_m\|_\infty \leq \|F\|_\infty$. The proof is completed by using the maximum modulus principle and \eqref{eq:bhpoly} with the weaker constant $e^m$. 
	\end{proof}
\end{lem}

We shall require the following number theoretic estimate, due to Balazard \cite{balazard1989remarques}. (Alternatively, we could have used the weaker Lemma 4.2$'$ of \cite{KQ01} which would have sufficed for our purpose.) 
\begin{lem}
	\label{lem:bala} Let $R(x,y,m) = \{n \leq x \,:\, P^-(n)>y,\, \Omega(n)=m\}$. For $x \geq y \geq 2$ there is some absolute positive constant $c$ such that 
	\begin{equation} \label{eq:rest}
		\Phi(x,y,M) = \sum_{m \geq M} \left|R(x,y,m)\right| \leq \frac{x}{y^M} \left(\log{x}\right)^{y} e^{cy}. 
	\end{equation}
	\begin{proof}
		This is a slightly weaker version of Corollary 1 in \cite{balazard1989remarques}. 
	\end{proof}
\end{lem}

\section{Proof of the upper bound} We may suppose $\|f\|_\infty=1$ without loss of generality. Let $x \geq y \geq 2$ and split $f(s)$ as in \eqref{eq:nsplit}, such that $\|\widehat{f}\|_1 = \sum_\kappa \|\widehat{f_\kappa}\|_1$. Furthermore, $\|f_\kappa\|_\infty\leq 1$ by Lemma \ref{lem:nsplit} and the indices $n$ of the sums $f_\kappa(s)$ are all $y$-rough numbers less than $x$. Suppose that we can prove, for some particular $y$ and each $y$-smooth $\kappa$, that 
\begin{equation}
	\label{eq:splitest} \|\widehat{f_\kappa}\|_1 \leq \sqrt{x} \exp\left(\left(-1+3\frac{\log\log\log{x}}{\log\log{x}}+\mathcal{O}\left(\frac{1}{\log\log{x}}\right)\right)\sqrt{\frac{\log{x}\log\log{x}}{2}}\right). 
\end{equation}
Combining \eqref{eq:splitest} with the splitting \eqref{eq:nsplit} we then obtain
\[\|\widehat{f}\|_1 \leq \sqrt{x} \exp\left(\left(-1+3\frac{\log\log\log{x}}{\log\log{x}}+\mathcal{O}\left(\frac{1}{\log\log{x}}\right)\right)\sqrt{\frac{\log{x}\log\log{x}}{2}}\right)\sum_{\kappa}1.\]
The sum is taken over all possible $y$-smooth $\kappa$ such that there is some $n\leq x$ with $n =\kappa\eta$, for some $y$-rough $\eta$. 

If $\kappa \leq x$ is $y$-smooth, its $\pi(y)$ possible prime factors have exponents between $0$ and $\Omega(\kappa) \leq \log{x}/\log{2}$. Hence we estimate\footnote{This crude estimate is sufficient, since the term $(\log{x})^y$ also appears when using \eqref{eq:rest}.}
\[\sum_{\kappa} 1  = \left|\left\{ \kappa \leq x \, : \, P^+(\kappa) \leq y \right\}\right|\leq \left(1+\frac{\log{x}}{\log{2}}\right)^{\pi(y)} \ll \left(\log{x}\right)^y = \exp(y\log\log{x}).\]
To force this contribution to be contained in the error term, we choose $y=\sqrt{\log{x}/(\log\log{x})^3}$. We now aim to prove \eqref{eq:splitest}, and may assume that the indices $n$ in the sum $f(s)$ are $y$-rough and that $\|f\|_\infty \leq 1$. Furthermore, since $x$ and $y$ now are fixed, we let $R_m=R(x,y,m)$.

Let us split the sum of the coefficients according to the number of prime factors of the index $n$. By combining Hölder's inequality with Lemma \ref{lem:homoineq} we obtain
\[\|f\|_1 = \sum_{m} \sum_{n \in R_m}|a_n| \leq \sum_{m} |R_m|^\frac{m-1}{2m} \left(\sum_{n \in R_m} |a_n|^\frac{2m}{m+1}\right)^\frac{m+1}{2m} \leq \sum_{m} |R_m|^\frac{m-1}{2m} e^m.\]
However, this straightforward approach has some problems.

In particular, when $y\log\log{x} \gg M\log{y}$, the estimate \eqref{eq:rest} of Lemma 4 is worse than the trivial bound $\Phi(x,y,M)\leq x$. We remedy this by letting $\alpha>0$ and
\[M_1 = \alpha\sqrt{\frac{\log{x}}{\log\log{x}}}.\]
For $m\leq M_1$ we will estimate using $|R_m|\leq x$. In particular
\[\Sigma_1 = \sum_{m \leq M_1} \sum_{n \in R_m} |a_n| \leq \sum_{m \leq M_1} x^\frac{m-1}{2m} e^m \leq \sqrt{x} \exp\left(-\frac{\log{x}}{2M_1}+ \mathcal{O}\left(\sqrt{\frac{\log{x}}{\log\log{x}}}\right)\right).\]
The main term is
\[-\frac{\log{x}}{2M_1} = - \frac{1}{2\alpha}\sqrt{\log{x}\log\log{x}}.\]
This implies that the largest value we may choose is $\alpha=1/\sqrt{2}$.

The use of Lemma \ref{lem:homoineq} comes at the cost of the factor $e^m$. To ensure that $e^m$ is contained in the error term, we will use the Cauchy--Schwarz inequality for all $m\geq M_2$, where $\beta>0$ and
\[M_2 = \beta \sqrt{\frac{\log{x}}{\log\log{x}}}.\]
We use the Cauchy--Schwarz inequality and extend the $\ell^2$-sum of the coefficients to obtain
\[\Sigma_3 = \sum_{m \geq M_2} \sum_{n \in R_m} |a_n| \leq \left(\sum_{m \geq M_2} \sum_{n \in R_m} 1\right)^\frac{1}{2}\left(\sum_{n\leq x} |a_n|^2\right)^\frac{1}{2} \leq \sqrt{\Phi(x,y,M_2)}.\]
The final inequality follows by the fact that $\|\widehat{f}\|_2 \leq \|f\|_\infty \leq 1$, which is evident by orthogonality in view of the Bohr lift. Furthermore, using Lemma \ref{lem:bala} we estimate
\[\sqrt{\Phi(x,y,M_2)} \leq \left(\frac{x}{y^{M_2}}(\log{x})^ye^{cy}\right)^\frac{1}{2} \leq \sqrt{x}\exp\left(-\frac{M_2\log{y}}{2}+\mathcal{O}\left(\sqrt{\frac{\log{x}}{\log\log{x}}}\right)\right).\]
Here, the main term is
\[-\frac{M_2\log{y}}{2} = -\frac{\beta}{4}\sqrt{\log{x}\log\log{x}} + \frac{3\beta}{4}\log\log\log{x}\sqrt{\frac{\log{x}}{\log\log{x}}}.\]
This implies that the smallest $\beta>0$ we may take is $\beta = 2 \sqrt{2}$.

What remains is to consider the values $M_1 < m < M_2$, which are of the form
\[m = \gamma \sqrt{\frac{\log{x}}{\log\log{x}}},\]
where $\alpha = 1/\sqrt{2} < \gamma < 2\sqrt{2} = \beta$. We follow the procedure of the first step, but for these $m$ the estimate $|R_m| \leq \Phi(x,y,m)$ is sharper than the trivial bound $|R_m| \leq x$. In particular, by Hölder's inequality we obtain
\[\Sigma_2 = \sum_{M_1 < m < M_2} \sum_{n \in R_m} |a_n| \leq \sum_{M_1 < m < M_2}|R_m|^\frac{m-1}{2m} e^m \leq \sum_{M_1 < m < M_2} \Phi(x,y,m)^\frac{m-1}{2m} e^m.\]
As in the first step, the number of summands and the factor $e^m$ are absorbed in the error term. Using Lemma \ref{lem:bala} we estimate
\[\Phi(x,y,m)^\frac{m-1}{2m} \leq \sqrt{x} \exp\left(-\frac{\log{x}}{2m}-\frac{m}{2}\log{y}+\mathcal{O}\left(\sqrt{\frac{\log{x}}{\log\log{x}}}\right)\right).\]
The main terms are
\[-\frac{\log{x}}{2m}-\frac{m}{2}\log{y} = - \left(\frac{1}{2\gamma}+\frac{\gamma}{4}\right)\sqrt{\log{x}\log\log{x}}+\frac{3\gamma}{4}\log\log\log{x}\sqrt{\frac{\log{x}}{\log\log{x}}}.\]
The first term is maximal at $\gamma = \sqrt{2}$, while the second term is maximal at $\gamma = 2\sqrt{2}$, yielding the required terms. We write 
\[\|\widehat{f}\|_1 = \sum_{m} \sum_{n \in R_m} |a_n| =  \Sigma_1+\Sigma_2+\Sigma_3\]
and using the estimates given above to complete the proof of \eqref{eq:splitest}.\qed

\section{Acknowledgments} 
The author would like to express his gratitude to K. Seip for suggesting the topic of the paper and helpful discussions along the way, and to H. Queffélec and S. Schlüters for pertinent remarks.
\bibliographystyle{amsplain} 
\bibliography{sidonest} 
\providecommand{\bysame}{\leavevmode\hbox to3em{\hrulefill}\thinspace}
\providecommand{\MR}{\relax\ifhmode\unskip\space\fi MR }
\providecommand{\MRhref}[2]{%
  \href{http://www.ams.org/mathscinet-getitem?mr=#1}{#2}
}
\providecommand{\href}[2]{#2}
\begin{thebibliography}{10}

\bibitem{balazard1989remarques}
Michel Balazard, \emph{Remarques sur un theoreme de {G. Halasz} et {A.
  Sarkozy}}, Bull. Soc. Math. France \textbf{117} (1989), no.~4, 389--413.

\bibitem{bpsbohr}
F.~Bayart, D.~Pellegrino, and J.~Seoane-Sep{ú}lveda, \emph{The {B}ohr radius
  of the $n$-dimensional polydisk is equivalent to $\sqrt{(\log{n})/n}$},
  arXiv:1310.2834 (2013).

\bibitem{BH31}
H.~F. Bohnenblust and Einar Hille, \emph{On the absolute convergence of
  {D}irichlet series}, Ann. of Math. \textbf{32} (1931), no.~3, 600--622.

\bibitem{bohr-darstellung}
Harald Bohr, \emph{Darstellung der gleichmäßigen konvergenzabszisse einer
  {Dirichletschen} reihe {$\sum_{n=1}^\infty \frac{a_n}{n^s}$} als funktion der
  koeffizienten der reihe}, Archiv der Mathematik und Physik \textbf{21}
  (1913), no.~3, 326--330.

\bibitem{B13b}
\bysame, \emph{Über die {B}edeutung der {P}otenzreihen unendlich vieler
  {V}ariabeln in der {T}heorie der {D}irichletschen {R}eihen}, Nachr. Akad.
  Wiss. Göttingen Math.-Phys. Kl. (1913), 441--488.

\bibitem{B13a}
\bysame, \emph{Über die gleichmässige {K}onvergenz {D}irichletscher
  {R}eihen}, J. Reine Angew. Math. \textbf{143} (1913), 203--211.

\bibitem{bombieri2009kahane}
Enrico Bombieri and Jean Bourgain, \emph{On {K}ahane's ultraflat polynomials},
  Journal of the European Mathematical Society \textbf{11} (2009), no.~3,
  627--703.

\bibitem{Br08}
R{é}gis de~la Bret{è}che, \emph{Sur l'ordre de grandeur des polynômes de
  {D}irichlet}, Acta Arith. \textbf{134} (2008), no.~2, 141--148.

\bibitem{bh-hyper}
Andreas Defant, Leonhard Frerick, Joaquim Ortega-Cerd{\`a}, Myriam
  Ouna{\"{\i}}es, and Kristian Seip, \emph{The {B}ohnenblust-{H}ille inequality
  for homogeneous polynomials is hypercontractive}, Ann. of Math. (2)
  \textbf{174} (2011), no.~1, 485--497.

\bibitem{erdos}
Paul Erd\H{o}s, \emph{An inequality for the maximum of trigonometric
  polynomials}, Ann. Polon. Math \textbf{12} (1962), 151--154.

\bibitem{kahane1980polynomes}
Jean-Pierre Kahane, \emph{Sur les polyn{\^o}mes {\'a} coefficients
  unimodulaires}, Bulletin of the London Mathematical Society \textbf{12}
  (1980), no.~5, 321--342.

\bibitem{KQ01}
S.~V. Konyagin and H.~Queff{é}lec, \emph{The translation {$\frac12$} in the
  theory of {D}irichlet series}, Real Anal. Exchange \textbf{27} (2001/02),
  no.~1, 155--175.

\bibitem{Queffelec}
H.~Queff{é}lec, \emph{H. {B}ohr's vision of ordinary {D}irichlet series; old
  and new results}, J. Anal. \textbf{3} (1995), 43--60.

\end{thebibliography}
\end{document}